\documentclass[11pt]{article}
\usepackage{amsthm,amsmath,amsfonts,latexsym,amssymb,mathrsfs,ulem}
\usepackage{color}

\usepackage{amscd, amsfonts, graphicx,tikz, color, environ}
\usepackage{amssymb}
\usepackage{amsmath}
\usepackage{amsfonts}
\newtheorem{lemma}{Lemma}[section]
\newtheorem{proposition}[lemma]{Proposition}
\newtheorem{theorem}[lemma]{Theorem}
\newtheorem{corollary}[lemma]{Corollary}

\theoremstyle{definition}
\newtheorem{remark}[lemma]{Remark}

\newtheorem{problem}[lemma]{Problem}
\newtheorem{example}[lemma]{Example}

\theoremstyle{plain}

\newcommand{\supp}{\operatorname{supp}}

\numberwithin{equation}{section}

\renewcommand{\thefootnote}{\fnsymbol{footnote}}

\definecolor{darkgreen}{rgb}{.1,.5,0}

\begin{document}

\title{ A Recursive Approach to the Matrix Moment Problem}

\author{Ra\'ul E. Curto, Abderrazzak Ech-charyfy, Kaissar Idrissi, and El Hassan Zerouali}

\date{}

\maketitle


\noindent {\bf Abstract.} \ In this paper, we study the truncated matrix moment problem in one variable through linearly recursive matrix extensions. \ We give necessary and sufficient conditions for a linearly recursive matrix extension of finite data to be a matrix moment sequence in the classical cases of Hamburger, Stieltjes, and Hausdorff moment problems.   

\setcounter{page}{1}


\renewcommand{\thefootnote}{}
\footnote{
\\
\vskip -.5cm \ 

\noindent \hskip -.4cm \textit{Mathematics Subject
Classification
(2020).} {Primary 44A60, 47A57; Secondary 47B20, 47B37}\\
\noindent
 \textit{Keywords.} Matrix moment problem, linearly recursive relations, sequences of symmetric matrices} 


\section{Introduction}

The real and complex one-dimensional moment problems are well-known in classical \linebreak analysis and have been widely studied by many mathematicians and engineers since the late 19th century. \ Stieltjes studied the problem on the positive real line \cite{stieltjes1894recherches}, Hamburger on the real line \cite{hamburger1920erweiterung}, and Hausdorff on a closed and bounded interval \cite{hausdorff1921summationsmethoden}. \ The key to resolving the full moment problem lies in the use of the Riesz-Haviland Theorem on one hand (\cite{riesz1923probleme} and \cite{Haviland}), and an appropriate characterization of positive polynomials on the other hand. \ The solution to the one-dimensional moment problem is expressed in terms of the positivity of certain associated Hankel matrices.

In the truncated case, a more constructive approach is employed, relying on techniques from finite-dimensional linear algebra in conjunction with extension methods. \ For truncated moment problems, the first-named author and L. Fialkow led a comprehensive investigation centered around the concept of flat extensions; see \cite{curto1996solution,curto1998flat,curto2000truncated} for further information. \ In the matrix-valued case, D. Kimsey and H. Woerdeman provided a solution based on commutativity conditions of certain matrices \cite{kimsey2013truncated}. \ Several other mathematicians have also tackled the truncated moment problem using various approaches, such as K. Schm\"udgen's operator approach \cite{schmudgen2017moment} and F.H. Vasilescu's dimensional stability approach \cite{vasilescu2012dimensional}, among others. 

Matrix moment problems have become increasingly popular and have been actively studied in recent years, owing to their applications, which attract the attention of a growing number of researchers. \ In \cite{dette2002note}, H. Dette and W.J. Studden used matrix moment problem techniques to delve into a matrix version of the Quotient-Difference algorithm in numerical analysis. \ They demonstrated that this algorithm can be employed to derive the coefficients of recurrence relations for matrix orthogonal polynomials on intervals such as $[0, 1]$ and $[0,+\infty)$, based on their moment generating functional. \ Their paper includes illustrative examples that extend the classical orthogonal polynomials defined on the real line. \ 

In \cite{dette2005note}, Dette and Studden investigated the maximum value and the determinant of canonical Hankel matrices associated with Hausdorff matrix moment sequences. \ The obtained results were then used in optimal design problems within linear models. \ This highlights a statistical application that involves maximizing matrix Hankel determinants constructed from matrix moments. 

Furthermore, J.B. Lasserre combined global optimization with polynomials and addressed the moment problem in a broader context \cite{lasserre2001global}. \ He emphasized the interactions between truncated moment problems and global optimization techniques using sum-of-squares polynomials and other related techniques.

In the present paper, we adopt the linearly recursive relations approach widely studied by M. Rachidi et al in the one-dimensional case, as shown in  \cite{taher2001recursive} and \cite{chidume2001solving}. \ Our paper is organized as follows. \ In Section 2, we introduce the matrix moment problem and give the general preparatory results needed for its resolution. \ We exhibit the description of positive matrix-valued polynomials on the classical intervals $I=[0,1], I= [0,+\infty)$  and $I= {\mathbb R}$. \ Section 3 is devoted to the study of linearly recursive matrix sequences.\ We show that a linearly recursive matrix sequence admits a representing matrix charge if and only if the associated minimal polynomial has all its roots distinct and real. \ 

In Section 4, we characterize the linearly recursive matrix moment sequences in terms of the positivity of an associated finite Hankel matrix. \ We give several computational results in the case of linearly recursive matrix sequences of order $2$ and $3$. 

\section{Preliminaries} In this paper, we use the following notation and symbols. \ $\mathbb{R}$, $\mathbb{N}:=\{1, 2, 3, \ldots\}$, and $\mathbb{N}_0:=\mathbb{N} \cup \{0\}$ will denote the sets of real numbers, positive integers, and nonnegative integers, respectively. \ The polynomial algebra $\mathbb{R}[X]$ consists of all polynomials with real coefficients. \ For $P\in \mathbb{R}[X]$, ${\mathcal Z}(P)$ denotes the zeros set of $P$. \ The algebra of $p \times p$ matrices with real entries is denoted by $\mathcal{M}_{p}(\mathbb{R})$, where $I_p$ is the identity matrix and $0_p$ is the null matrix in $\mathcal{M}_{p}(\mathbb{R})$. \ We use $\mathcal{\mathcal S} _{p}(\mathbb{R})$ to denote the subspace of symmetric matrices in $\mathcal{M}_{p}(\mathbb{R})$ and  $\mathcal{\mathcal S} _{p}^+(\mathbb{R})$ to denote the set of nonnegative matrices in $\mathcal{\mathcal S} _{p}(\mathbb{R})$, equipped with the Loewner order, denoted by "$\succeq$". \ We will write $\mathcal{M}_{p}(\mathbb{R}[X])$  for the algebra of $p \times p$ matrices with  entries in $\mathbb{R}[X$]. \ Alternatively, $\mathcal{M}_{p}(\mathbb{R}[X])$ can be viewed as $\mathcal{M}_{p}(\mathbb{R})[X]$, the algebra of polynomials with  coefficients in $\mathcal{M}_{p}(\mathbb{R})$. \ In other words, every polynomial $P \in \mathcal{M}_{p}(\mathbb{R}[X])$ may take one of two different forms:
\[
P(X) = \sum_{k=0}^{n}A_{k}X^{k} = (P_{i,j})_{1\le i,j\le p},
\]
where $A_{k}=(a_{ij}^{(k)})_{1\le i,j\le p}\in \mathcal{M}_{p}(\mathbb{R})$ for $0\le k\le n$, and 
\[
P_{i,j}(X)=\displaystyle\sum_{k=0}^{n}a_{ij}^{(k)}X^{k} \in \mathbb{R}[X] 
\]
for every $1 \le i,j\le p$. \ If $A_n \ne 0_p$, we will say that $P$ is of degree $n$ and write $\deg (P)=n$. \\

Before introducing the matrix $K$--moment problem and the tracial $K$--moment problem, as stated in \cite{le2019tracial}, we recall that $\mathcal{B}(\mathbb{R})$ denotes the $\sigma$--algebra of Borel subsets in $\mathbb{R}$. \ A matrix-valued function $\mu : \mathcal{B}(\mathbb{R}) \to \mathcal{M}_p(\mathbb{R})$ is called a matrix charge if  
\[
\mu(B) := (\mu_{ij}(B))_{i,j=1,\ldots,p} \in \mathcal{M}_p(\mathbb{R}) \quad (\mbox{ for every } B \in \mathcal{B}(\mathbb{R})),
\]
where $\mu_{ij}$ is a scalar charge on $\mathcal{B}(\mathbb{R})$, for $1\le i, j \le p$.\\

A matrix charge  $\mu$ is called a matrix measure if for all $B \in \mathcal{B}(\mathbb{R})$ and for all $v \in \mathbb{R}^p$, we have $v^T \mu(B) v \geq 0$; or, equivalently, if $\mu(B) \in \mathcal{\mathcal S} _{p}^+(\mathbb{R})$. 

The support of a matrix charge  $\mu= (\mu_{ij})_{i,j=1,\ldots,p}$, denoted by $\supp(\mu)$, is defined as the smallest closed subset $F \subseteq \mathbb{R}$ such that $\mu(S)=0_{p}$ for every Borel set $S \subseteq \mathbb{R} \setminus F$. \ In the case of a matrix measure, the support coincides with the smallest closed subset $F \subseteq \mathbb{R}$ such that $\mu(\mathbb{R} \setminus F)=0_{p}$.

For a given matrix charge $\mu \equiv (\mu_{ij})_{i,j=1,\ldots,p}$, it is not difficult to prove the following general characterization of the support:
\[
\supp(\mu) = \bigcup_{i,j=1}^p \supp(\mu_{ij}).
\]
Moreover, if $\mu$ is a  matrix measure, then for every $B \in \mathcal{B}(\mathbb{R})$ and every $i,j=1,\dots,p$, we have
$
\begin{pmatrix}
    \mu_{ii}(B) & \mu_{ij}(B)\\
    \mu_{ij}(B) & \mu_{jj}(B)
\end{pmatrix} \succeq 0.
$
It follows, in particular, that $\mu_{ij}(B)^2 \leq \mu_{ii}(B)\mu_{jj}(B)$. \ This implies that $$\supp(\mu_{ij})\subseteq \supp(\mu_{ii})\cap \supp(\mu_{jj}),$$ \ and then 
\[
\supp(\mu) = \bigcup_{i=1}^p \supp(\mu_{ii}).
\]

An important class of matrix charges consists of those with finite support, called finitely atomic matrix charges. \ They are written as
\begin{equation}\label{fam}
\mu = \sum_{i=1}^{k} T_i \delta_{\lambda_i},
\end{equation}
where $\lambda_1, \ldots, \lambda_k \in \mathbb{R}$ and $T_1, \ldots, T_k \in \mathcal{M}_{p}(\mathbb{R})$.

Let $\mu=(\mu_{ij})_{1\le i,j\le p}$ be a matrix charge. \ A function $f : \mathbb{R} \to \mathbb{R}$ is called $\mu$--measurable if $f$ is $\mu_{ij}$--measurable for every $i, j = 1, \ldots, p$. \ In this case, the matrix-valued integral of $f$ with respect to the matrix charge $\mu$ is defined by
\[
\int_\mathbb{R} f(x) d\mu(x) := \left(\int_\mathbb{R} f(x) d\mu_{ij}(x)\right)_{i,j=1,\ldots,p} \in \mathcal{M}_p(\mathbb{R}).
\]
In particular, if $\mu$ is a finitely atomic matrix charge of the form \eqref{fam}, then 
\[
\int_\mathbb{R} f(x) d\mu(x) = \sum_{i=1}^{k}  f(\lambda_i)T_i\in \mathcal{M}_p(\mathbb{R}).
\]

Let $ {\mathcal S} =(S_{n})_{n\in \mathbb{N}}$ be a sequence in $\mathcal{\mathcal S} _{p}(\mathbb{R})$. \ We consider the real-valued linear form $\mathbf{L}_{{\mathcal S}}$ on $\mathcal{M}_{p}(\mathbb{R}[X])$ defined by
\[\begin{array}{rcll}
   \mathbf{L}_{{\mathcal S}} : &\mathcal{M}_{p}(\mathbb{R}[X]) &\longmapsto &  \mathbb{R} \\
     & P = \sum_{k=0}^{m}A_{k}X^{k} &\mapsto &\mathbf{L}_{{\mathcal S}}(P) := \sum_{k=0}^{m}\operatorname{tr}\left(A_{k}S_{k}\right),
\end{array}
\]
where $\operatorname{tr}(A)$ is the trace of the matrix $A$.

The linear form $\mathbf{L}_{{\mathcal S}}$ is known as the Riesz functional associated with the matrix sequence ${\mathcal S}$. \ Given a closed set $K\subseteq \mathbb{R}$, the tracial $K$--moment problem is  introduced as follows.
\begin{problem}\label{trmmp}
The tracial $K$--moment problem asks for necessary and sufficient conditions for the existence of a  matrix measure $\mu=(\mu_{ij})_{1\leq i,j \leq p}$ satisfying the following conditions:
\begin{equation*}
    \supp(\mu)\subseteq K \quad \text{and} \quad \mathbf{L}_{{\mathcal S}}(P)=\int_{K}\text{tr}\big(P(x)d\mu(x)\big) \text{ for every } P\in\mathcal{M}_{p}(\mathbb{R}[X]).
\end{equation*}
\end{problem}
The matrix $K$--moment problem is formulated as follows.

\begin{problem}\label{mmp} 
The matrix $K$--moment problem asks for necessary and sufficient conditions for the existence of a matrix measure $\mu=(\mu_{ij})_{1\leq i,j \leq p}$ satisfying
\begin{equation*}
    \supp(\mu)\subseteq K \quad \text{and} \quad S_{k}=\int_{K}x^{k}d\mu(x) \text{ for every } k\in\mathbb{N}.
\end{equation*}
\end{problem}
We retain the usual definitions in the case of scalar moment problems. \ A matrix sequence (finite matrix sequence) for which the corresponding  $K$--matrix moment problem (truncated matrix moment problem) has a solution is referred to as a  $K$--matrix moment sequence (resp. truncated $K$--matrix moment sequence). \ Problems \ref{trmmp} and \ref{mmp} are referred to as the \textit{Hamburger matrix moment problems} if $K=\mathbb{R}$, the \textit{Stieltjes matrix moment problems} if $K=[0,+\infty)$, and the \textit{Hausdorff matrix moment problems} if $K=[0,1]$.
   
Let us now recall a useful result from \cite{le2019tracial}.
\begin{theorem}[$\text{\cite[Proposition 1 and Theorem 2.1]{le2019tracial}}$]\label{x}
Let $\mathcal{S} \equiv (S_\alpha)_{\alpha \in \mathbb{N}_0}$ be a sequence of symmetric matrices in $\mathcal{M}_p(\mathbb{R})$ and let $K$ be a closed set in $\mathbb{R}^n$. \ The following statements are equivalent:
\begin{enumerate}
     \item ${\mathcal S}$ is a   $K$--matrix moment sequence;
    \item $\mathbf{L}_{{\mathcal S}}$ is a tracial $K$--moment functional;
    \item $\mathbf{L}_{{\mathcal S}}(P)\geq 0$ for all $P\in\mathcal{P}(K)$, where $\mathcal{P}(K) = \{P \in \mathcal{M}_p(\mathbb{R}[X]) \mid P_{|K} \succeq 0\}$.
\end{enumerate}
\end{theorem}
\begin{remark}
   In the truncated case, Theorem \ref{x}  has been extended by C. M{\"a}dler and K. Schm{\"u}dgen to all measure spaces \cite{madler2023truncated}.
\end{remark}

It follows from the previous discussion that a necessary condition for ${\mathcal S} $ to be a $K$--matrix moment sequence can be expressed as follows:
\begin{equation}\label{ham+}
\mathbf{L}_{{\mathcal S}}(PP^t) \geq 0, \text{ for every } P\in \mathcal{M}_p(\mathbb{R}[X]);
\end{equation}
or, equivalently, 
\begin{equation}\label{ham}
    \sum_{i,j=0}^{n}\operatorname{tr}(A_iS_{i+j}A_{j}^t) = \operatorname{tr}\left(\sum_{i,j=0}^{n}A_iS_{i+j}A_{j}^t\right) \geq 0,
\end{equation}
for every $A_0,\dots ,A_n\in\mathcal{M}_p(\mathbb{R})$.

\begin{remark}\label{1a}
   Condition \eqref{ham+} is equivalent to  $H_{m}: = (S_{i+j})_{0\leq i,j\leq m}\succeq 0$,  for every $m\in \mathbb{N}$. \ Thus, for an arbitrary set $K$ of real numbers, a necessary solubility condition for the $K$--matrix moment problem is $H_{m}\succeq 0$ for every $m\in\mathbb{N}$.
\end{remark}
As stated above, (\ref{ham+}) is a necessary condition for $\mathcal S$ to be a $K$--matrix moment sequence. \ In the case of $K = \mathbb{R}$, (\ref{ham+}) is also a sufficient condition. \ For more specific results, we have the following theorem.  
\begin{theorem}$($\cite[Corollary 1]{cimprivc2013moment}$).$\label{solution} Let ${\mathcal S}=(S_{n})_{n\in \mathbb{N}}$ be a sequence in $\mathbb{R}$.  For every $k, m \in \mathbb{N}$, we denote by $H_{m}^{(k)}=(S_{i+j+k})_{0\leq i,j\leq m}$ the moment block matrix, and we let $H_{m}:=H_{m}^{(0)}$. \ Then:
\begin{enumerate}
    \item ${\mathcal S}$ is a Hamburger matrix moment sequence if and only if $H_{m}\succeq 0$ for every $m\in\mathbb{N}$ (in which case $K = \mathbb{R}$).
    \item ${\mathcal S}$ is a Stieltjes matrix moment sequence if and only if $H_{m}\succeq 0$ and $H_{m}^{(1)}\succeq 0$ for every $m\in\mathbb{N}$ (in which case $K=[0,+\infty))$.
    \item ${\mathcal S}$ is a Hausdorff matrix moment sequence if and only if $H_{m}\succeq 0$, $H_{m}^{(1)}\succeq 0$, $H_{m}-H_{m}^{(1)}\succeq 0$, and $H_{m}^{(1)}- H_{m}^{(2)}\succeq 0$ for every $m\in \mathbb{N}$ (in which case $K = [0,1]$).
\end{enumerate}
\end{theorem}
 Recall that  $P\in\mathcal{M}_{p}\left(\mathbb{R}[X]\right)$ is said to be a sum-of-squares matrix polynomial if there exists $P_1,\cdots, P_k\in\mathcal{M}_{p}\left(\mathbb{R}[X]\right)$ such that $P(X) = P_1(X)P_{1}^{t}(X) +\cdots + P_k(X)P_{k}^{t}(X)$, where $P^{t}$ is the transpose of the matrix polynomial $P$. \ The most common approach to obtain a proof of Theorem \ref{solution} is to combine Theorem \ref{x} with the following description of positive matrix polynomials.
\begin{proposition}[$\text{\cite[Proposition 3]{cimprivc2013moment}}$]\label{n}
For $P\in\mathcal{M}_{p}\left(\mathbb{R}[X]\right)$, the following statements hold.
\begin{enumerate}
     \item $P\succeq 0$ on $\mathbb{R}$ if and only if $P$ is a sum of two square matrix polynomials.
     \item $P\succeq 0$ on $\left[0,+\infty \right[ \iff P=A+XB$, where $A$ and $B$ are sum-of-squares matrix polynomials.
     \item $P\succeq 0$ on $\left[0,1 \right]$ if and only if $P=A+XB+(1-X)C+X(1-X)D$, where $A, B, C, \mbox{ and } D$ are sum-of-squares matrix polynomials.
\end{enumerate}
\end{proposition}
Another solution to the Hausdorff matrix moment problem is obtained by combining Theorem \ref{solution} above together with the description of positive matrix polynomials on $[0,1]$ and by using the Bernstein decomposition, as provided in \cite[Corollary 3.7]{le2019tracial}. 
\begin{proposition}\label{sh}
    A matrix sequence $(S_n)_{n\ge 0}$ is a Hausdorff matrix moment sequence if and only if 
    $\sum\limits_{i=1}^{k} (-1)^i\binom{k}{i}\text{tr}(S_{i+l}G)\ge 0,$
    for all positive definite matrices $G$ and $k,l \in \mathbb{N}.$
\end{proposition}

Thanks to the expressions $X=X^2+X(1-X)$ and $1-X=(1-X)^2+X(1-X)$, we derive the following finer characterization of non-negative matrix polynomials, improving assertion 3 of Proposition \ref{n}.
\begin{proposition}\label{ab}
Let $P\in\mathcal{M}_{p}\left(\mathbb{R}[X]\right)$ \ The following statements are equivalent: \  

(i) \ $P\succeq 0 \ on \ [0,1]$;

(ii) \ $P(X)=A+X(1-X)B$, where $A$ and $B$ are sum-of-squares matrix polynomials.

\end{proposition}
As a consequence of these results, we obtain the following solution of the Hausdorff matrix moment problem, with a simpler condition than in statement 3 of Theorem \ref{solution}.
\begin{theorem}
The following statements are equivalent:
\begin{enumerate}
    \item ${\mathcal S}$ is a Hausdorff matrix moment sequence;
    \item $H_{m}\succeq 0$ and $H_{m}^{(1)}-H_{m}^{(2)}\succeq 0$ for every $m\in \mathbb{N}$;
    \item $ \sum\limits_{i=1}^{k} (-1)^i\binom{k}{i}H_{m}^{(i)}\succeq 0,$  for every $m\in \mathbb{N}.$  
\end{enumerate}
\end{theorem}

\section[Moment Problems and Recursive Matrix Sequences]{The Truncated Matrix Moment Problem and Linearly Recursive Matrix Sequences}\label{second}
In the truncated matrix moment problem, only finite data \, ${\mathcal S} = (S_0,\cdots, S_{r-1})$ of symmetric matrices is provided. \
The truncated tracial $K$--moment problem and the truncated matrix $K$--moment problem are introduced as follows.
\begin{problem}
The truncated tracial $K$--moment problem asks for necessary and sufficient conditions for the existence of a  matrix measure $\mu$ 
with  $\supp(\mu)\subseteq K$ and such that
\begin{equation*}
    \mathbf{L}_{{\mathcal S}}(P)=\int_{K}\text{tr}\big(P(x)d\mu(x)\big), \mbox{ for all } P\in \mathcal{M}_{p}(\mathbb{R}[X]) \mbox{ with } \ \deg(P)\le r-1.
\end{equation*}
\end{problem}
\begin{problem}\label{tmmp} 
The truncated matrix $K$--moment problem asks for necessary and sufficient conditions for the existence of a matrix measure $\mu$ with  $\supp(\mu)\subseteq K $ and such that
\begin{equation*}
    \quad S_{k}=\int_{K}x^{k}d\mu(x) \ \ for \ \ 0\le k\le r-1.
\end{equation*}
\end{problem}
The Richter–Tchakaloff Theorem (cf. \cite[Thm. 1.24]{schmudgen2017moment}) asserts that each scalar moment functional on a finite-dimensional vector space of functions admits a finitely atomic representing measure. \ The Richter–Tchakaloff Theorem has been extended recently to the matrix case by C. M{\"a}dler and K. Schm{\"u}dgen in \cite[Theorem 5.1]{madler2023truncated}. \ It follows that the above-mentioned truncated matrix moment problems, the full matrix moment problem, and the matrix moment problems for linearly recursive matrix sequences are all equivalent. \ This last fact motivates our restriction, in the remainder of this paper, to sequences of matrices given by $r$ initial data items and satisfying a linear recurrence relation.

It turns out that, under an additional condition on the recurrence relation of order $r$, all the solubility conditions of the classical matrix moment problem can be reduced to a unique condition, expressed as the positivity of $H_{r-1}$. 

To be more precise, let $a_{0},a_{1},\ldots,a_{r-1}$ be some fixed real numbers with $a_{r-1} \neq 0$ and consider the sequence of matrices ${\mathcal S}=(S_{n})_{n\in\mathbb{N}}$ in $\mathcal{M}_{p}(\mathbb{R})$ defined by the following linear recurrence relation:
\begin{equation}\label{xx}
    S_{n+1}=a_{0} S_{n}+a_{1} S_{n-1}+\ldots+a_{r-1} S_{n-r+1},  \quad n\geq r-1,
\end{equation}
where   $r\ge 1$ is an integer and $ S_{0}, S_{1},\ldots, S_{r-1} \in \mathcal{M}_{p}(\mathbb{R})$ are the initial conditions. \ In this paper, we will refer to sequences satisfying \eqref{xx} as linearly recursive matrix sequences.\\

The scalar case, $p=1$, has been widely studied in the literature. \ These sequences are called $r$--generalized Fibonacci sequences and have been a key point in solving several known problems; see, for example,  \cite{taher2001recursive,dubeau1994weighted} and \cite{WahbiRachidi}.

For the matrix case, we  write ${\mathcal S} =(S_n)_{n=1}^\infty=((s_{ij}^{(n)})_{1\leq i,j \leq p})_{n=1}^\infty$. \ Then for $n\in \mathbb{N}$, it follows  that \eqref{xx} is equivalent to the  next family of linear scalar relations:
\begin{equation}\label{re}
   s_{ij}^{(n+1)}=a_{0}s_{ij}^{(n)}+a_{1} s_{ij}^{(n-1)}+\dots+a_{r-1}s_{ij}^{(n-r+1)},
\end{equation}
for every $  n\geq r-1,\text{ and  } 1\leq i,j \leq p.$

The polynomial $Q_{\mathcal S}(X):=X^{r}-a_{0}X^{r-1} -\dots-a_{r-2}X-a_{r-1}$ is called a characteristic polynomial associated with the linearly recursive matrix sequence ${\mathcal S}$.

As shown in the following example, a linearly recursive matrix sequence always has many characteristic polynomials. 
\begin{example}\label{ex1}
Let $M\in \mathcal{M}_{p}(\mathbb{R})$, and $S_{n}:=(n+1)M$. \ Then the sequence ${\mathcal S} =(S_{n})_{n\in\mathbb{N}}$ verifies the following two linearly recursive matrix relations.
\[
\begin{array}{lll}
    S_{0}=M, \; \; S_{1}=2M, \; \; S_{2}=3M \;\;  &\text{ and } \; S_{n+1}=S_{n}+S_{n-1}-S_{n-2}, &\text{ for } n\ge 2,  \\
    S_{0}=M, \; \; S_{1}=2M  \; &\text{ and } \;   S_{n+1}=2S_{n}-S_{n-1}, &\text{ for } n\ge 1.
\end{array}
\]
It follows that $P(X)=X^{3}-X^{2}-X+1$ and $Q(X)= X^{2}-2X+1$ are two characteristic polynomials for the sequence ${\mathcal S} =(S_{n})_{n\in\mathbb{N}}$.
\end{example}

Let $\mathcal{P}_{\mathcal S} $ be the vector space generated by all characteristic polynomials associated with the linearly recursive matrix sequence ${\mathcal S} =(S_{n})_{n\in\mathbb{N}}$. \ We have the following result. 
\begin{proposition}\label{idea}
For every linearly recursive matrix sequence ${\mathcal S} =(S_{n})_{n\in\mathbb{N}}$, the set $\mathcal{P}_{\mathcal S} $ is an ideal in $\mathbb{R}[X]$. \ In particular, there exists a unique monic characteristic polynomial $P_{\mathcal S} \in\mathcal{P}_{\mathcal S} $ with minimal degree $r$ such that $\mathcal{P}_{\mathcal S} =P_{\mathcal S} \mathbb{R}[X]$. 
\end{proposition}
\begin{proof}
 
Let $\mathcal{M}_{p}(\mathbb{R})^{\mathbb{N}}$ be the vector space of all sequences $s \equiv (M_n)_{n\in\mathbb{N}}$ with values in $\mathcal{M}_{p}(\mathbb{R})$ and let ${\mathcal L}(\mathcal{M}_{p}(\mathbb{R})^{\mathbb{N}})$ be the algebra of linear maps on $\mathcal{M}_{p}(\mathbb{R})^{\mathbb{N}}$. \ The shift operator $W$   on  $\mathcal{M}_{p}(\mathbb{R})^{\mathbb{N}}$ is  defined by
\begin{equation}\label{endo}
    W((M_{n})_{n\in\mathbb{N}}):=(M_{n+1})_{n\in\mathbb{N}}.
\end{equation} 
Consider also the algebra homomorphism
$$
\begin{array}{ccc}
   \Phi: \mathbb{R}[X] & \to & {\mathcal L}(\mathcal{M}_{p}(\mathbb{R})^{\mathbb{N}})\\
     Q =\displaystyle\sum_{i=0}^{m}a_{i}X^{i}  &\mapsto &   Q(W)= \displaystyle\sum_{i=0}^{m}a_{i}W^{i}.
\end{array}
$$
Using \eqref{xx}, observe that $Q \in \ker(\Phi)$ if and only if $(0_p)_{n \in \mathbb{N}}=Q(W)(\mathcal S)$, and if and only if $Q \in \mathcal{P}_{\mathcal S}$. 

Clearly, for $Q_1, Q_2 \in \ker(\Phi)$ and $P\in \mathbb{R}[X]$, we have $Q_1+ Q_2 \in \ker(\Phi)$ and $PQ_1 \in \ker(\Phi)$. \ It follows that $ \mathcal{P}_{\mathcal S} =\ker(\Phi)$ is an ideal in the principal ideal domain $\mathbb{R}[X]$. \ Hence, there exists a unique monic polynomial $P_{\mathcal S} \in \mathcal{P}_{\mathcal S} $ such that $\mathcal{P}_{\mathcal S} =P_{\mathcal S} \mathbb{R}[X]$.
\end{proof}

The monic polynomial $P_{\mathcal S} $ will be called the minimal polynomial associated with ${\mathcal S} =(S_{n})_{n\in \mathbb{N}}$. \ We will say in the sequel that a linearly recursive matrix sequence is of order $r$ if the associated minimal polynomial has degree $r$.

\subsection{Minimal polynomial and finitely atomic representing matrix charge}
We next investigate the relationship between the minimal polynomial of a linearly recursive matrix sequence and the minimal polynomials of its scalar entries. \ More precisely, we have the following:
\begin{theorem}
Let ${\mathcal S} =(S_n)_{n\in\mathbb{N}}=\big((s_{ij}^{(n)})_{1\leq i,j \leq p}\big)_{n\in\mathbb{N}}$ be a linearly recursive sequence with minimal polynomial $P_{\mathcal S} $. \ For $1\leq i,j\leq p$, let $P_{ij}$ be the minimal polynomial associated with the linearly recursive scalar sequence  $s_{ij}=(s_{ij}^{(n)})_{n\in\mathbb{N}}$. \ Then
\begin{equation*}
    P_{\mathcal S} =lcm\big\{P_{ij}\ , \ 1\leq i,j\leq p\big\},
\end{equation*}
where $lcm\big\{P_{ij} \ , \ 1\leq i,j\leq p\big\}$ stands for the least common multiple of the family of polynomials $\big\{P_{ij} \ ,\   1\leq i,j\leq p\big\}$.
\end{theorem}
\begin{proof}
Let $P_{lcm}$ be the least common multiple of the family of polynomials $\big\{P_{ij}  \ , \ 1\leq i,j\leq p\big\}$, and let us show that $P_{\mathcal S} =  P_{lcm}$.
\ Since $P_{\mathcal S} $ is a characteristic polynomial for $s_{ij}=(s_{ij}^{(n)})_{n\in\mathbb{N}}$, it is a multiple of $P_{ij}$ for every $1\leq i,j\leq p$, and then $P_{\mathcal S} $ is a multiple of $lcm\big\{P_{ij}  \ , \  1\leq i,j\leq p\big\} = P_{lcm}$. \ Conversely, $P_{lcm}$ is a multiple of $P_{ij}$ for every $1\leq i,j\leq p$ and then is a characteristic polynomial of $s_{ij}=(s_{ij}^{(n)})_{n\in\mathbb{N}}$ for every $1\leq i,j\leq p$. \ Now, since addition and scalar multiplication act linearly on entries, it follows that $P_{lcm}$ is a characteristic polynomial for ${\mathcal S}$. \ We deduce that $P_{lcm}$ is a multiple of $P_{\mathcal S} $, which completes the proof.
\end{proof}
\begin{remark}
  For an alternative proof of the previous theorem, one can observe that for every $1\leq i,j\leq p$, the polynomial $P_{ij}$ generates the ideal  $P_{ij}\mathbb{R}[X]$, and that 
  $$
  P_{\mathcal S} \mathbb{R}[X]=\displaystyle\bigcap_{i,j=1}^{p} P_{ij}\mathbb{R}[X] = P_{lcm}\mathbb{R}[X].
  $$
  \end{remark}

We derive some additional useful properties.
\begin{corollary}\label{aaa}Using the above notation, the following statements hold:
\begin{enumerate}
    \item  ${\mathcal Z}(P_{\mathcal S} )=\displaystyle\bigcup_{ i,j=1}^{ p}{\mathcal Z}(P_{ij})$;
    \item $P_{\mathcal S} $ has distinct roots if and only if $P_{ij}$ has distinct roots for every $1\le i,j\le p.$
\end{enumerate}
\end{corollary}
Now, we present a cornerstone result in solving the matrix moment problem associated with linearly recursive matrix sequences. \ This result characterizes the existence of representing charges for linearly recursive matrix sequences. 
\begin{theorem}\label{fla}
Let ${\mathcal S} =(S_{n})_{n\in\mathbb{N}}$ be a linearly recursive matrix sequence with a minimal polynomial $P_{\mathcal S} $. \ The following statements are equivalent.
\begin{enumerate}
    \item ${\mathcal S}$ admits a representing matrix charge $\mu$;
        \item $P_{\mathcal S} $ has real distinct roots. 
\end{enumerate}  
Moreover, when a representing charge $\mu$ exists, we have $\supp(\mu)={\mathcal Z}(P_{\mathcal S})$.
\end{theorem}
\begin{proof}
Let ${\mathcal S} =(S_{n})_{n\in\mathbb{N}}$ be a  matrix sequence. \ It is clear that ${\mathcal S} =(S_{n})_{n\in\mathbb{N}}$  is linearly recursive if and only if the scalar  sequences $s_{ij}=(s_{ij}^{(n)})_{n\in\mathbb{N}}$ are linearly recursive for all $1\leq i, j\leq p$. 

    Similarly, the  matrix sequence ${\mathcal S} =(S_{n})_{n\in\mathbb{N}}$ admits a representing matrix charge $\mu=(\mu_{ij})_{1\leq i,j\leq 1}$ if and only if the scalar  sequence $s_{ij}=(s_{ij}^{(n)})_{n\in\mathbb{N}}$ admits a representing scalar charge $\mu_{ij}$ for every $1\leq i, j\leq p$. 
    
    Hence if ${\mathcal S} =(S_{n})_{n\in\mathbb{N}}$ is linearly recursive (and so are $s_{ij}=(s_{ij}^{(n)})_{n\in\mathbb{N}}$), according to \cite[Proposition 2.4]{taher2001recursive}, this is equivalent to 
    the fact that $P_{ij}$ has distinct roots for every $1\leq i, j\leq p$. \ By  Corollary  \ref{aaa}, the latter is equivalent to    $P_{\mathcal S} $ having distinct roots.
    
    Moreover, $\supp(\mu_{ij})={\mathcal Z}(P_{ij})$ for every $1\leq i, j\leq p$, and then $$\supp(\mu)=\displaystyle\bigcup_{ i,j=1}^{ p}\supp(\mu_{ij})=\bigcup_{ i,j=1}^{ p}{\mathcal Z}(P_{ij})={\mathcal Z}(P_{\mathcal S} ).$$
\end{proof}

\begin{example}\begin{enumerate} \item  For $n \in \mathbb{N}$, let  $S_n:=\begin{pmatrix}
1&n\\
n&1
\end{pmatrix}$, and let ${\mathcal S} :=(S_n)_{n\in\mathbb{N}}$. \ It is straightforward to verify that $S_{n+2}=2 S_{n+1}-S_n \; \; (n \in \mathbb{N})$, and that the   minimal polynomial is $P_{\mathcal S} (X):=(X-1)^2$, which has no distinct roots. \ Thus, the sequence ${\mathcal S} $ admits no representing matrix charge.
\item Let ${\mathcal S} =(S_n)_{n\ge 0}$ be a matrix sequence  such that $S_{n+2}=-S_n$. \ Then $P_{\mathcal S} (X)= X^2+1$  and since ${\mathbb R}\cap{\mathcal Z}(P_{\mathcal S} )=\emptyset$, we have   $S_n$ admits no representing matrix charge.
\item Denote  $A=\begin{pmatrix}
0&1\\
1&0
\end{pmatrix}$, and  $S_{n}=A^{n}$ for $n\in \mathbb{N}$. \ The characteristic polynomial of $A$ is given by $P_{A}(X):=X^{2}-1=(X+1)(X-1)$. \ Using the Cayley-Hamilton Theorem, we get $A^n(A^2-I_2)=0_2$  for every $n\ge 0$ and hence
\begin{equation*}
    S_{n}=S_{n-2}, \quad n\geq 2 \text{ with } S_{0}=I_{2} \text{ and } S_{1}=A.
\end{equation*}
In particular, $P_A$ is the minimal polynomial associated with $(S_n)_{n\in\mathbb{N}}$, and then $(S_n)_{n\in\mathbb{N}}$ admits a representing matrix charge. \  We obtain easily
\[
S_{n}=\frac{1}{2}\begin{pmatrix}
1-(-1)^{n}&1+(-1)^{n}\\
1+(-1)^{n}&1-(-1)^{n}
\end{pmatrix} = 
\frac{1}{2}\begin{pmatrix}
 1&1\\
 1&1
\end{pmatrix}+
\frac{(-1)^{n}}{2}\begin{pmatrix}
 -1&1\\
 1&-1
\end{pmatrix}.
\]
We can express $S_{n}$ as
\[
S_{n}=\int_{\supp(\mu)}t^{n}d\mu(t),
\]
where
\[
\mu=\begin{pmatrix}
\mu_{11}&\mu_{12}\\
\mu_{21}&\mu_{22}
\end{pmatrix}= 
\frac{1}{2}\begin{pmatrix}
 1&1\\
 1&1
\end{pmatrix}\delta_{1}+
\frac{1}{2}\begin{pmatrix}
 -1&1\\
 1&-1
\end{pmatrix}\delta_{-1},
\]
and  $\supp(\mu)={\mathcal Z}(P_{A})=\{-1,1\}$. \\
Here,   $\mu$ is  not a matrix measure since $ \begin{pmatrix}
 -1&1\\
 1&-1
\end{pmatrix}$ is not positive.
\end{enumerate}
We mention at this stage that the Cayley-Hamilton Theorem shows that $A^n$ is a linearly recursive sequence for every matrix $A\in M_p(\mathbb{R})$ associated with the characteristic polynomial   $P_A(X)= {\rm det}(A-XI_p)$. \ For a symmetric matrix $A$, its minimal polynomial has only real simple roots, and $A$ is similar to a diagonal matrix in $M_p(\mathbb{R})$. \ It follows that a  representing matrix charge always exists for $(A^n)_{n\ge 0}$, with $A$ a symmetric matrix.
\end{example}

\subsection{Minimal polynomial and finitely atomic representing matrix measure}
In the remainder of this paper, we will focus on linearly recursive matrix sequences ${\mathcal S} = (S_{n})_{n\in \mathbb{N}}$ whose minimal polynomial $P_{\mathcal S} $ has all its roots distinct and real. \ We will write ${\mathcal Z}(P_{\mathcal S} )=\{\lambda_{1},\lambda_{2},\ldots,\lambda_{r}\}$ and denote $\mu = \sum_{q=1}^{r}T_{q}\delta_{\lambda_{q}}$, with $T_{1},\dots, T_{r} \in \mathcal{M}_{p}(\mathbb{R})$ an associated finitely atomic matrix charge provided by  Theorem \ref{fla}.

Our next objective is to determine necessary and sufficient conditions on the initial data of ${\mathcal S} =(S_0,\cdots, S_{r-1})$ to guarantee that $\mu$ is a finitely atomic matrix measure. \ In other words, when are the coefficients $T_{1},\dots, T_{r} $  in $\mathcal{\mathcal S} _{p}^+(\mathbb{R})$? \ According to Remark \ref{1a}, this occurs precisely when the Hankel matrix $H_n$ satisfies the condition
   \begin{equation}\label{cnf}
       H_n=(S_{i+j})_{0\le i,j\le n}\succeq 0 \quad (\mbox{ for every } n\in\mathbb{N}).
   \end{equation}
As in the scalar case \cite[Proposition 3.2]{taher2001recursive}, the following proposition provides a simple equivalent condition to $\eqref{cnf}$ in the case of linearly recursive matrix sequences of order $r$.
\begin{proposition}\label{prop2}
Let ${\mathcal S} $ be a linearly recursive matrix sequence of order $r$, and suppose that its minimal polynomial $P_{\mathcal S} $ has only real distinct roots. \ The following statements are equivalent:
\begin{enumerate}
    \item $H_{r-1} \succeq 0$;
    \item $H_n \succeq  0$ for every $n \in\mathbb{N}$.
\end{enumerate}
\end{proposition}
\begin{proof} We only need to show $(1)\Rightarrow (2)$. \  To this aim,
denote by ${\mathcal S} =(S_n)_{n\in\mathbb{N}}$ the linearly recursive matrix sequence of order $r$. \ Let ${\mathcal Z}(P_{\mathcal S} ) = \{\lambda_1, \lambda_2, \ldots, \lambda_r\}$ be the set of distinct roots of $P_{\mathcal S} $, and let $\mu$ be an associated finitely atomic matrix charge as in Theorem \ref{fla}. \ We have
\begin{equation*}
    \mu = \sum_{q=1}^{r} T_{q}\delta_{\lambda_{q}} \, \text{ for some } \, T_{1},\dots, T_{r} \in \mathcal{M}_{p}(\mathbb{R}).
\end{equation*}
Assume that $H_{r-1} \succeq 0$ and let $\vec{c_0}, \vec{c_1}, \ldots,\vec{c_n} \in \mathbb{R}^{p}$. \ We write
\begin{equation}\label{po}
\sum_{i,j=0}^{n}\vec{c_i} \, S_{i+j} \, \vec{c_j}^{T} = \sum_{q=1}^{r} \vec{P}(\lambda_q) T_q \vec{P}(\lambda_q)^{T},
\end{equation}
where $\vec{P}(X) = \sum_{k=0}^{n} \vec{c}_k X^k\in \mathbb{R}^p_n[X]$ denotes a vector polynomial with coefficients in $\mathbb{R}^{p}$.  \ By setting  $\vec{c}_k = (c_{k,i})_{1\leq i\leq p}$ for $0\leq k\leq n$, we obtain  $\vec{P}(X) = (P_i(X))_{1\leq i\leq p}$ with $P_n(X) = \sum_{k=0}^{n} c_{k,i} X^k \in \mathbb{R}_n[X]$. \
Now, for fixed $1\leq i\leq p$, there exist two polynomials $Q_i$ and $R_i$ such that
\begin{equation*}
P_i = P_{\mathcal S}Q_i + R_i, \text{ with } \text{deg}(R_i) < \text{deg}(P_{\mathcal S} ) = r. 
\end{equation*}
Hence, for every vector-valued polynomial, we have
\begin{equation*}
\vec{P} = P_{\mathcal S}.\vec{Q} + \vec{R},
\end{equation*}
with $\vec{Q}(X) = (Q_i(X))_{1\leq i\leq p}$ and $\vec{R}(X) = (R_i(X))_{1\leq i\leq p}= \sum_{k=0}^{r-1} \vec{\alpha}_k X^k$. \ Here  $\vec{\alpha}_0,   \ldots, \vec{\alpha}_{r-1} \in \mathbb{R}^{p}$ are the coefficients of the vector-valued polynomial $\vec{R}$. \ Moreover, for every $1\leq q\leq r$, we have  $\vec{P}(\lambda_q) = \vec{R}(\lambda_q)$, since $P_{\mathcal S} (\lambda_q) = 0$. \ Thus, according to the hypothesis and (\ref{po}), we have
\begin{equation*}
\sum_{i,j=0}^{n}\vec{c_i} \, S_{i+j} \, \vec{c_j}^T = \sum_{q=1}^{r} \vec{R}(\lambda_q) T_q \vec{R}(\lambda_q)^{T} = \sum_{i,j=0}^{r-1} \vec{\alpha_i} \, S_{i+j} \, \vec{\alpha_i}^{T} \geq 0.
\end{equation*}
It readily follows that $H_n \succeq 0$ for every $n \ge 0$.
\end{proof}
\begin{remark}
    Note that condition (2) in the previous proposition is equivalent to the positivity of the following infinite block moment matrix
\[H_\infty:= (S_{i+j})_{i, j\geq 0} = \left( \begin{matrix}
                S_0 & S_1 & S_2 & \ldots \\
                S_1 & S_2 & S_3 &\ldots \\
                S_2 & S_3 & S_4 &\ldots \\
                \vdots & \vdots & \vdots & \ddots
              \end{matrix}
              \right).\]
\end{remark}
We identify the polynomial $P_{n}(X) = \sum_{k=0}^{n}a_kX^k$ with the column matrix vector $\widehat{P}_{n} = (a_0I_{p}, \dots, a_nI_{p}, 0_{p}, 0_{p}, \dots)^T$, and we define the canonical inner product on $\mathbb{R}[X]$ by 
$$
\langle\sum_{k=0}^{n}a_kX^k, \sum_{k=0}^{m}b_kX^k\rangle = \sum_{k=0}^{min(n,m)}a_kb_k.
$$ 
Setting $H_\infty P := H_\infty \widehat{P}$, we can consider $H_\infty$ as a linear map on $\mathbb{R}[X]$. \ With this notation, we derive the following trivial properties.

\begin{lemma}\label{car}
Let ${\mathcal S} = (S_n)_{n \in \mathbb{N}}$ be a linearly recursive matrix sequence with associated minimal polynomial $P_{\mathcal S} $. \ Then we have:
\begin{enumerate}
    \item $P \in \mathcal{P}_S$ if and only if $H_\infty {P} = {0}$. \ In particular, we have $H_\infty {P}_{\mathcal S} = {0}$.
    \item For every $P, Q \in \mathbb{R}[X] $ we have $\langle Q,H_\infty {P} \rangle=\langle 1, H_\infty {PQ}\rangle.$
\end{enumerate}
\end{lemma}

\begin{proposition}\label{pri}Let ${\mathcal S} = (S_n)_{n \in \mathbb{N}}$ be a linearly recursive matrix sequence of order $r$ with associated minimal polynomial $P_{\mathcal S} $ and let $H_{r-1} := (S_{i+j})_{0 \leq i, j \leq r-1} \succeq 0$.  
\begin{enumerate}
\item The following statements are equivalent:
  \begin{enumerate}
  \item $P \in \mathcal{P}_{\mathcal S}$;
  \item $P^n \in \mathcal{P}_{\mathcal S}$ for every $n \in \mathbb{N}^*$;
  \item $P^m \in \mathcal{P}_{\mathcal S}$ for some $m \in \mathbb{N}^*$.
  \end{enumerate}
  \item The minimal polynomial $P_{\mathcal S} $ has distinct roots.
\end{enumerate}
\end{proposition}
\begin{proof}
$(1)$ Suppose that $H_{r-1} \succeq 0$. \ By Proposition \ref{prop2}, $H_n \succeq 0$ for every $n\ge 1$. \ Due to the fact that $\mathcal{P}_{\mathcal S}$ is an ideal of $\mathbb{R}[X]$, we have $(a) \Rightarrow (b)$, and it is evident that $(b) \Rightarrow (c)$ holds. \ Then, we only need to prove that $(c) \Rightarrow (a)$. \ 
Assume that there exists a positive integer $m$ such that $P^m \in \mathcal{P}_{\mathcal S}$. \ According to Lemma \ref{car}, this is equivalent to $H_\infty {P^m} = {0}$. \
If $m = 1$, then it is straightforward that $P \in \mathcal{P}_{\mathcal S}$. 

Now, for the case where $m \ge 2$, we use Lemma \ref{car} to get
\[
\langle P^{m-1}, H_\infty P^{m-1}\rangle = \langle{P^{m-2}}, H_\infty {P^{m}}\rangle = 0,
\]
 and since $H_\infty \succeq 0$ we obtain  $H_\infty {P^{m-1}} = 0$. \ By arguing recursively  on $m$, we deduce that $H_\infty {P} = {0}$, and therefore, that  $P \in \mathcal{P}_{\mathcal S}$.
\\$(2)$  Let $P_{\mathcal S} (X) = \displaystyle\prod_{i=1}^{k} (X -\lambda_i)^{m_i}$, and let $m\ge \displaystyle\max\{m_i \mid i=1, \ldots, k\}$. \ We have 
 $\displaystyle\left(\prod_{i=1}^{k} (X -\lambda_i)\right)^m \in \mathcal{P}_S$, since $\mathcal{P}_{\mathcal S}$ is an ideal of $\mathbb{R}[X]$. \ Using (1), we obtain $\displaystyle\prod_{i=1}^{k} (X -\lambda_i) \in \mathcal{P}_{\mathcal S}$, and then $P_{\mathcal S} (X) = \displaystyle\prod_{i=1}^{k} (X -\lambda_i)$, since $P_{\mathcal S} $ is the minimal polynomial.
 \end{proof}
\section{Linearly Recursive Moment Sequences}
In this section, we discuss the classical linearly recursive matrix moment sequences on ${\mathbb R}$.
\
\subsection{ Hamburger recursive matrix moment sequences}
 \begin{theorem}\label{Hamb}
Let ${\mathcal S} = (S_{n})_{n\in\mathbb{N}}$ be a linearly recursive matrix sequence of order $r$ with  $P_{\mathcal S} $ the associated minimal polynomial has real distinct roots. \ The following statements are equivalent:
\begin{enumerate}
    \item  ${\mathcal S} $ is a Hamburger matrix moment sequence;
    \item $H_{r-1} \succeq 0$.
\end{enumerate}
\end{theorem}
 \begin{proof}Let ${\mathcal S} = (S_{n})_{n\in\mathbb{N}}$ be a linearly recursive matrix sequence of order $r$, and let $P_{\mathcal S} $ be its minimal polynomial. \ The direct implication is clear. \ Now, we assume that $H_{r-1}\succeq 0$. \ Since the minimal polynomial $P_{\mathcal S} $ has $r$ real  distinct roots,  ${\mathcal Z}(P_{\mathcal S} ) = \{\lambda_1, \dots, \lambda_r\}$, from Proposition \ref{prop2} we derive that ${\mathcal S} = (S_{n})_{n\in\mathbb{N}}$ admits a finitely atomic generating charge $\mu$, such that $\mu = \displaystyle\sum_{i=1}^{r}T_{i}\delta_{\lambda_{i}}$ for some $T_1, \dots, T_{r} \in M_{p}(\mathbb{R})$. \ Thus, it remains to show that $T_i \in \mathcal{\mathcal S} _{p}^+(\mathbb{R})$ for $i = 1, \ldots, r$. \ To this aim, since $H_{r-1}\succeq 0$, we have:
 \begin{equation}\label{pos}
    \displaystyle\sum_{i,j=0}^{r-1}\vec{c}_{i} \, S_{i+j} \, \vec{c}_{j}^T = \sum_{i=0}^{r}\vec{P}(\lambda_{i})T_{i}\vec{P}(\lambda_{i})^T \ge 0,
\end{equation}
 for every vector polynomial $\vec{P}(X) = \displaystyle\sum_{k=0}^{r-1}\vec{c}_{k}X^{k}$, where $\vec{c}_{k} \in \mathbb{R}^{p}$.
 
Now, we fix $i = 1, \ldots, r$, and we consider the interpolating Lagrange polynomial
 \begin{equation*}
    L_{i}(X) := \prod_{\substack{j=1\\j\neq i}}^{r}\frac{(X-\lambda_{j})}{(\lambda_{i}-\lambda_{j})}.
\end{equation*}
 Also, for an arbitrary vector $\vec{v}$ in $\mathbb{R}^p$, we define $\vec{Q}_i(X): = \vec{v}\cdot L_i(X)$. \ It is clear that $\vec{Q}_i(\lambda_j) = \vec{v}\delta_{ij}$, where $\delta_{ij}$ is the Kronecker delta. \ Applying \eqref{pos} to the vector polynomial $\vec{P} = \vec{Q}_i$, we obtain
 $  \vec{v}^T \, T_{i} \, \vec{v} \ge 0, $ and then $T_i \in \mathcal{\mathcal S} _{p}^+(\mathbb{R})$ for every $1\le i \le r$. 
 \end{proof}
 \begin{remark}
  \ We can recover the previous theorem as a consequence of Theorem \ref{solution} and Proposition \ref{prop2}.
 \end{remark}
  As a corollary, we have the next completion result.
\begin{corollary}\label{mr} Let $r\ge 1$ and let $ {\mathcal S} \equiv (S_0, S_1, \cdots, S_{r-1})$ be a finite family of symmetric matrices. \ The following statements are equivalent:
\begin{enumerate}
     \item ${\mathcal S} $ is a Hamburger truncated moment sequence;
     \item There exists a family $(S_{r}, \cdots, S_{ 2r-2 })$ of \, $r-1$ symmetric matrices such that \linebreak 
     $(S_{i+j})_{0\le i, j \le r-1} \succeq 0.$
\end{enumerate}   
\end{corollary}
\subsection{Stieltjes and  Hausdorff recursive matrix sequences}
Using similar proofs, we obtain the next result.
\begin{theorem}\label{stie}
Let ${\mathcal S} = (S_{n})_{n\in\mathbb{N}}$ be a linearly recursive symmetric matrix sequence of order $r$ with  $P_{\mathcal S} $ the associated minimal polynomial has real distinct roots. \  Then
\begin{enumerate}
    \item  ${\mathcal S} $ is a Stieltjes matrix moment sequence if and only if $H_{r-1} \succeq 0$ and $H^{(1)}_{r-1} \succeq 0$.
    \item  ${\mathcal S} $ is a Hausdorff matrix moment sequence if and only if $H_{r-1} \succeq 0$ and $H^{(1)}_{r-1}-H^{(2)}_{r-1} \succeq 0$.
\end{enumerate}
\end{theorem}
We further deduce the following corollary.
\begin{corollary} Let $r\ge 1$ and $ {\mathcal S}=(S_0, S_1, \cdots, S_{r-1})$ be a finite family of symmetric matrices. \ The following statements are equivalent:
\begin{enumerate}
     \item ${\mathcal S} $ is a Stieltjes truncated matrix moment sequence;
     \item There exists a family $(S_{r}, \cdots, S_{ 2r-1 })$ of \, $r$ symmetric matrices such that $(S_{i+j})_{0\le i, j \le r-1} \succeq 0$ and $ (S_{i+j+1})_{0\le i, j \le r-1} \succeq 0.$
\end{enumerate}   
\end{corollary}
\begin{corollary} Let $r\ge 1$ and $ {\mathcal S} = (S_0, S_1, \cdots, S_{r-1})$ be a finite family of symmetric matrices. \ The following statements are equivalent:
\begin{enumerate}
     \item $ {\mathcal S} $ is a Hausdorff truncated matrix moment sequence;
     \item There exists $(S_{r}, \cdots, S_{ 2r})$ a family of symmetric matrices such that $(S_{i+j})_{0\le i, j \le r-1} \succeq 0$ and $ (S_{i+j+1}- S_{i+j+2})_{0\le i, j \le r-1}\succeq 0.$
\end{enumerate}   
\end{corollary}
\begin{remark}   Theorem \ref{Hamb} and Theorem \ref{stie} have been addressed in \cite{kimsey2013truncated} by  D. Kimsey and H. Woerdeman, who obtain necessary and sufficient conditions for the existence of minimal finitely atomic representing measures (with the proper conditions on their support), while also handling the case of odd degree. \ In the present paper, we adopt a recursive approach that yields a simple solution for both the odd and even cases in one variable. \ A suitable adaptation of the present results to the case of several variables, which we do not pursue here, is likely to provide an extension of the results given in \cite{kimsey2013truncated}.
\end{remark}

\subsection{Some special cases}
Next, we investigate two cases of linearly recursive matrix sequences of special interest. \ Without any loss of generality, we will take $S_{0}=I_{p}.$
\subsubsection{Case r=2.}  
Let ${\mathcal S} =(S_{n})_{n\in \mathbb{N}}$ be a linearly recursive matrix sequence, with initial data $(I_p,S_1)$, such that $S_{n+1}=a_{0}S_{n}+a_{1}S_{n-1}$ for $n\ge 1$, where $a_{0}^2+4a_{1}>0$. \ We write the associated minimal polynomial as $P_{\mathcal S} (X)=(X-\lambda_{0})(X-\lambda_{1})$ with $\lambda_{0}<\lambda_{1}$. \ According to Theorem \ref{fla}, the sequence ${\mathcal S} $ admits a finitely atomic generating matrix charge $\mu= T_{0}\delta_{\lambda_{0}}+T_{1}\delta_{\lambda_{1}}$, for some $T_{0},T_{1} \in \mathcal{M}_{p}(\mathbb{R})$. \ In particular, we have
\[
\left\{
\begin{array}{lllcl}
T_{0} + T_{1} &=& I_{p}, \\
\lambda_{0}T_{0} + \lambda_{1}T_{1} &=& S_{1},
\end{array}
\right.
\]
and then
\begin{equation}
T_0 = \frac{1}{\lambda_1-\lambda_0}(\lambda_1I_p-S_1) \quad \text{ and } T_1 = \frac{1}{\lambda_1-\lambda_0}(S_1-\lambda_0I_p). \label{T0T1}
\end{equation}

We now derive the following result.

\begin{theorem}\label{posit} \ Let $\lambda_0$ and $\lambda_1$ be real numbers such that $\lambda_0 < \lambda_1$, let $S_1 \in \mathcal{M}_p(\mathbb{R})$, and let $\sigma := T_0\delta_{\lambda_0}+T_1 \delta_{\lambda_1}$, where $T_0, T_1 \in \mathcal{M}_p(\mathbb{R})$ are as in  (\ref{T0T1}). \ The following statements are equivalent:
\begin{enumerate}
    \item $\mu$ is a matrix measure;
    \item $\lambda_0I_2\preceq S_{1}\preceq \lambda_1I_2$.
\end{enumerate}
\end{theorem}

The next corollary highlights the impact of the choice of  $\lambda_0$ and $\lambda_1$.
\begin{corollary}\label{coor} Let $A\in \mathcal{\mathcal S} _{p}(\mathbb{R})$. \ Then: 
\begin{enumerate}
    \item $(I_{p}, A)$ is always  a truncated Hamburger moment problem;
    \item $(I_{p}, A)$ is a truncated Stieltjes moment problem if and only if $A\in \mathcal{\mathcal S} ^+_{p}(\mathbb{R})$;
    \item $(I_{p}, A)$ is a truncated Hausdorff moment problem if and only if $A$ is   a contraction in $ \mathcal{\mathcal S} ^+_{p}(\mathbb{R})$.
\end{enumerate}
\end{corollary}

\begin{proof}
For (1), using Theorem \ref{posit}, it suffices to choose $\lambda_0= -\|A\|$  and $  \lambda_1=\|A\|$ so that $\lambda_0I_p \preceq A \preceq \lambda_1I_p$ and define ${\mathcal S} =(S_n)_{n\in\mathbb{N}}$ as follows:
\begin{equation*}
    S_{n+1} = (\lambda_0 + \lambda_1)S_n - \lambda_0\lambda_1S_{n-1}, \quad n\ge 1,
\end{equation*}
with initial conditions $S_0=I_p$ and $S_1=A$.\\
We obtain the second assertion by remarking that necessarily $\lambda_0 \in {\mathbb R}^+$, and the third by observing that we should have $\lambda_0, \lambda_1 \in [0,1]$.
\end{proof}

\begin{remark}We notice here that for every $A\in\mathcal{\mathcal S} _{p}(\mathbb{R})$, the sequence ${\mathcal S} =(A^n)_{n\in\mathbb{N}}$ is a linearly recursive matrix sequence that admits a finitely atomic representing matrix charge $\mu_A$. \ However, $\mu_A$ is a matrix measure if and only if $A\in\mathcal{\mathcal S} ^+_{p}(\mathbb{R})$. \ Thus, the previous corollary deals with a more general setting. 
\end{remark}

In the scalar case $p=1$, we recover Proposition 4.1 of \cite{taher2001recursive} in the following corollary.

\begin{corollary}
Let ${\mathcal S} =(s_{n})_{n\in\mathbb{N}}$ be a linearly recursive real sequence such that $s_{0}=1$, $s_{1}\in \mathbb{R}$, and $s_{n+1}=a_{0}s_{n}+a_{1}s_{n-1}$, with $a_0, a_1 \in \mathbb{R}$. \ Also, let $P_{\mathcal S} (X):=X^2-a_0X-a_1$. \ The following statements are equivalent: 
\begin{enumerate}
    \item ${\mathcal S} $ is a real moment sequence;
    \item $P_{\mathcal S} (s_1)\leq 0$.
\end{enumerate}

\end{corollary}

 \subsubsection{Case  r=3.}\label{432} \ Let ${\mathcal S} = (I_p,S_1,S_2)$ be a given initial data of symmetric matrices. \ Using  Corollary \ref{mr}, we have that ${\mathcal S} $ is a truncated moment sequence if and only if there exist $S_3$ and  $S_4$ such that $H_2\succeq 0$. \ It is then necessary to have
 $S_2\succeq S_1^2$.  Thus, as in the scalar case, there exist initial data $(I,S_1,S_2)$ without any matrix moment sequence extension. \ See \cite{taher2001recursive}, for the scalar case. \ We exhibit below a constructive method to obtain a representing matrix measure when $(I,S_1,S_2)$ is a truncated moment sequence.
 
Let ${\mathcal S} = (S_n)_{n \in \mathbb{N}}$ be a linearly recursive matrix sequence of order $r = 3$, with initial data $S_0 = I_p$, and $S_1, S_2 \in \mathcal{\mathcal S} _p(\mathbb{R})$. \ Let $P_{\mathcal S} $ be the associated minimal polynomial, such that $P_{\mathcal S} (X) = (X - \lambda_0)(X - \lambda_1)(X - \lambda_2)$, with $\lambda_0 < \lambda_1 < \lambda_2$. \ Denote  $\mu =  T_0\delta_{\lambda_0} +  T_1\delta_{\lambda_1} +  T_2 \delta_{\lambda_2}$for the associated finitely atomic representing matrix charge obtained from Theorem \ref{fla}. \ The coefficients $T_1, T_2, T_3 \in \mathcal{M}_p(\mathbb{R})$ will satisfy the following Vandermonde system:
$$
\left \{
\begin{array}{lllcl}
T_0& + T_1& + T_2 &=& I_p,\\
\lambda_0 T_0& + \lambda_1 T_1& + \lambda_2 T_2 &=& S_1,\\
\lambda_0^2 T_0& + \lambda_1^2 T_1& + \lambda_2^2 T_2 &=& S_2.
\end{array}
\right.
$$
Solving this linear matrix system, we obtain the expression of $T_0, T_1,$ and $T_2$   
$$\left \{\begin{array}{lll} T_0&=&(\lambda_3-\lambda_2){\Delta}.(S_{2}-S_{1}(\lambda_1+\lambda_2)+\lambda_1\lambda_2I_p),\\
T_1&=&(\lambda_2-\lambda_0){\Delta}.(-S_{2}+S_{1}(\lambda_0+\lambda_2)-\lambda_0\lambda_2I_p),\\
T_2&=&(\lambda_1-\lambda_0){\Delta}.(S_{2}-S_{1}(\lambda_0+\lambda_1)+\lambda_0\lambda_1I_p),
\end{array}
\right.
$$
where 
$$
\Delta=\dfrac{1}{(\lambda_2-\lambda_1)(\lambda_2-\lambda_0)(\lambda_1 -\lambda_0)}>0.
$$
Hence, we derive the following result.

\begin{proposition}\label{fin} \ With the notation as in Subsection \ref{432} above, the following statements are equivalent:
\begin{enumerate}
    \item $\mu$ is matrix measure;
\item $\left \{
\begin{array}{l}
S_{2}-(\lambda_1+\lambda_2)S_{1}+\lambda_1\lambda_2I_p\succeq 0, \\
-S_{2}+(\lambda_0+\lambda_2)S_{1}-\lambda_0\lambda_2I_p\succeq 0,  \\
S_{2}-(\lambda_0+\lambda_1)S_{1}+\lambda_0\lambda_1I_p \succeq 0; 
\end{array}
\right.$
\item $\left \{
\begin{array}{l}
\lambda_1(S_{1}-\lambda_2I_p)\le S_2-\lambda_2S_1\le \lambda_0(S_1-\lambda_2I_p), \\
\lambda_1(S_{1}-\lambda_0I_p)\le S_2-\lambda_0S_1\le \lambda_2(S_1-\lambda_0I_p).
\end{array}
\right.$
\end{enumerate}
\end{proposition}
 
\bigskip
Notice that  some necessary conditions to determine if $\mu$ is a matrix measure are given by
$$
S_2\succeq  S_1^2,  \ \ \lambda_0I_p\preceq S_{1}\preceq \lambda_2I_p  \mbox{ and } \lambda_{0}^2I_p\preceq S_{2}\preceq \lambda_{2}^2I_p.
$$ 
In the extremal cases $S_1= \lambda_0I_p$ or $S_1= \lambda_2I_p$, $(I_p, S_1, S_2)$ is a truncated matrix moment sequence if and only if  $ S_2= S_1^2$. \ This corresponds to the one-atomic matrix measure  $S_n= A \lambda^n$, for some $A\in\mathcal{\mathcal S} ^+_p(\mathbb{R}) $ and $\lambda \in {\mathbb R}$ as in the scalar case.
\begin{example}\label{ex3}
Let ${\mathcal S} =(S_{n})_{n\in \mathbb{N}}$ be a linearly recursive sequence of $\mathcal{\mathcal S} _{3}^{+}(\mathbb{R})$  satisfying $S_{n+1}=6S_{n}-10S_{n-1}+4S_{n-2}$ for $n\ge 2$, and with  initial data  $$S_{0}= \begin{pmatrix}
1 & 0 & 0  \\
0 & 1 & 0  \\
0 & 0 & 1 
\end{pmatrix}, \ S_1=\begin{pmatrix}
2 & -1 & 0  \\
-1 & 2 & -1  \\
0 & -1 & 2 
\end{pmatrix}\ \mbox{ and } \ S_2=\begin{pmatrix}
5 & -4 & 1  \\
-4 & 6 & -4  \\
1 & -4 & 5 
\end{pmatrix}.$$   
The associated minimal polynomial is given by
\begin{equation*}
    P_{\mathcal S} (X)=X^3-6X^2+10X-4=(X-\lambda_0)(X-\lambda_1)(X-\lambda_2),
\end{equation*}
where $\lambda_{0}:=2-\sqrt{2}$, $\lambda_{1}:=2$ and $\lambda_{2}:=2+\sqrt{2}$. \ We have
\begin{equation*}
\left\{\begin{array}{ll}
      S_{2}-(\lambda_1+\lambda_2)S_{1}+\lambda_1\lambda_2I_3  & =\begin{pmatrix}
1 & \sqrt{2} & 1  \\
\sqrt{2} & 2 & \sqrt{2}  \\
1 & \sqrt{2}& 1 
\end{pmatrix} \succeq 0, \\
   -S_{2}+(\lambda_0+\lambda_2)S_{1}-\lambda_0\lambda_2I_3  &  =\begin{pmatrix}
1 & 0 & -1  \\
0 & 0 & 0  \\
-1 & 0& 1 
\end{pmatrix} \succeq 0,\\
   S_{2}-(\lambda_0+\lambda_1)S_{1}+\lambda_0\lambda_1I_3 & =\begin{pmatrix}
1 & -\sqrt{2} & 1  \\
-\sqrt{2}  & 2 & -\sqrt{2}   \\
1 & -\sqrt{2} & 1 
\end{pmatrix} \succeq 0.
\end{array} \right.
\end{equation*}
Therefore, according to Proposition \ref{fin}, the sequence ${\mathcal S}= (S_n)_{n \in \mathbb{N}}$ admits  finitely atomic representing matrix measure $\mu$. \ Moreover, 
\begin{equation*}  \mu=T_0\delta_{\lambda_0}+T_1\delta_{\lambda_1}+T_2\delta_{\lambda_2},
\end{equation*}
where $$T_0 =\displaystyle\begin{pmatrix}
\frac{1}{4} & \frac{\sqrt{2}}{4} & \frac{1}{4}  \\
 \frac{\sqrt{2}}{4}  &\frac{1}{2}& \frac{\sqrt{2}}{4}  \\
\frac{1}{4} & \frac{\sqrt{2}}{4}& \frac{1}{4} 
\end{pmatrix}, \ T_1 =\displaystyle\begin{pmatrix}
\frac{1}{2} & 0 & -\frac{1}{2}  \\
0 & 0 & 0  \\
-\frac{1}{2} & 0 & \frac{1}{2}
\end{pmatrix}  \mbox{   and } T_2 =\displaystyle\begin{pmatrix}
\frac{1}{4} & -\frac{\sqrt{2}}{4}& \frac{1}{4}  \\
-\frac{\sqrt{2}}{4} & 2& -\frac{\sqrt{2}}{4}  \\
\frac{1}{4} & -\frac{\sqrt{2}}{4}& \frac{1}{4} 
\end{pmatrix}.$$ 
In particular, we have $S_n = T_0\lambda_0^n + T_1\lambda_1^n + T_2\lambda_2^n$ for every $n \in \mathbb{N}$.

\end{example}
\begin{remark}
If ${\mathcal S} =(S_n)_{n\in \mathbb{N}}=\big((s_{i,j}^{(n)})_{1\leq i,j \leq p}\big)_{n\in\mathbb{N}}$ is a matrix moment sequence, then   $(s_{ii}^{(n)})_{n\in\mathbb{N}}$ is also a scalar moment sequence for every $i=1,\dots,p$. \
Moreover, Example \ref{ex3} shows that this is no longer true for $s_{ij}$ when $i\ne j$. \ In the reverse direction, the sequence 
$$
S_n:= \begin{pmatrix}
1+2^n & 1+3 \cdot 2^n  \\
1+3 \cdot 2^n  & 2^n   \\
\end{pmatrix}=  \begin{pmatrix}
1 & 1    \\
1  & 0   \\
\end{pmatrix}1^n+ \begin{pmatrix}
1 &3    \\
3  &1  \\
\end{pmatrix}2^n
$$
has as entries real moment sequences. \ Since the associated matrix charge 
$$
\mu=  \begin{pmatrix}
1 & 1    \\
1  & 0   \\
\end{pmatrix}\delta_1+\begin{pmatrix}
1 &3    \\
3  &1  \\
\end{pmatrix}\delta_2
$$
is not a matrix measure, ${\mathcal S} $ is not a matrix moment sequence.
\end{remark}

\section{Declarations}
\subsection{Funding} \ The first-named author was partially supported by NSF grant DMS-2247167.  \ The last-named author was partially supported by the Arab Fund Foundation Fellowship Program. The Distinguished Scholar Award - File 1026. \ He also acknowledges the Mathematics Department at The University of Iowa for its kind hospitality during the preparation of this paper.

\subsection{Conflicts of interest/competing interests}

\subsection*{Data availability.}
All data generated or analyzed during this study are included in this article.

\vskip 1cm

Ra\'ul E. Curto

Department of Mathematics, University of Iowa, Iowa City, IA 52242,
U.S.A.

E-mail: raul-curto@uiowa.edu

\bigskip

Abderrazzak Ech-charyfy

Faculty of Sciences, Mohammed V University in Rabat, BP 1014 Rabat, Morocco.

E-mail: abderrazzak\_echcharyfy@um5.ac.ma
 
\bigskip

Kaissar Idrissi

Faculty of Sciences,  Mohammed V University in Rabat, BP 1014 Rabat, Morocco.

E-mail: kaissar.idrissi@fsr.um5.ac.ma

\bigskip

El Hassan Zerouali

Faculty of Sciences, Mohammed V University in Rabat,  BP 1014 Rabat, Morocco

E-mail: elhassan.zerouali@fsr.um5.ac.ma \& ezerouali@uiowa.edu 

\end{document}